\documentclass[aop,MSNbibl,nameyear,dvips]{arximspdf}
\usepackage{graphics}
\doi{10.1214/09-AOP508}
\volume{38}
\issue{3}
\pubyear{2010}
\firstpage{1263}
\lastpage{1279}

\makeatletter
\def\eqref#1{(\ref{#1})}
\newtheorem{theorem}{Theorem}
\newtheorem{lemma}[theorem]{Lemma}
\newtheorem{proposition}[theorem]{Proposition}
\newproclaim{remark}[theorem]{Remark}
\newcommand{\RRVlong}{Ram{\'{\i}}rez, Rider and Vir\'{a}g (\citeyear{RRV})}
\makeatother

\begin{document}
\begin{frontmatter}

\title{Large gaps between random eigenvalues}
\runtitle{Large gaps between random eigenvalues}

\begin{aug}
\author[A]{\fnms{Benedek} \snm{Valk\'{o}}\thanksref{a1}\ead[label=e1]{valko@math.wisc.edu}} \and
\author[B]{\fnms{B\'{a}lint} \snm{Vir\'{a}g}\thanksref{a2}\ead[label=e2]{balint@math.toronto.edu}\corref{}}
\runauthor{B. Valk\'{o} and B. Vir\'{a}g}
\thankstext{a1}{Supported in part by the Hungarian Scientific Research
Fund Grant K60708 and the NSF Grant DMS-09-05820.}
\thankstext{a2}{Supported by the Sloan and Connaught grants, the NSERC
discovery grant program and the Canada Research Chair program.}
\affiliation{University of Wisconsin--Madison and University of Toronto}
\address[A]{Department of Mathematics\\
University of Wisconsin--Madison\\ 480 Lincoln Dr\\Madison,
Wisconsin 53706\\
USA\\\printead{e1}} 
\address[B]{Departments of Mathematics\\
\quad and Statistics\\
University of Toronto\\
40 St. George Street\\
Toronto, Ontario M5S 2E4\\
Canada\\\printead{e2}}
\pdfauthor{Benedek Valko, Balint Virag}
\end{aug}

\received{\smonth{12} \syear{2008}}
\revised{\smonth{10} \syear{2009}}

%
\begin{abstract}
We show that in the point process limit of the bulk
eigenvalues of \mbox{$\beta$-ensembles} of random matrices, the
probability of having no eigenvalue in a fixed interval of
size $\lambda$ is given by
\[
\bigl(\kappa_\beta+o(1)\bigr)\lambda^{\gamma_\beta}\exp\biggl(-
\frac{\beta}{64}\lambda^2+\biggl(\frac{\beta}{8}-\frac14\biggr
)\lambda\biggr)
\]
as $\lambda\to\infty$, where
\[
\gamma_\beta=\frac{1}{4}\biggl(\frac{\beta}{2}+\frac{2}{\beta}-
3\biggr)
\]
and $\kappa_\beta$ is an undetermined positive constant.
This is a slightly corrected version of a prediction by Dyson
[\textit{J. Math. Phys.} \textbf{3} (1962) 157--165]. Our proof uses the new Brownian carousel
representation of the limit process, as well as the
Cameron--Martin--Girsanov transformation in stochastic
calculus.
\end{abstract}

%

\end{frontmatter}
%

\section{Introduction}

In the 1950s, Wigner endeavored to set up a probabilistic
model for the repulsion between energy levels in large
atomic nuclei. His first models were random meromorphic
functions related to random Schr\"{o}dinger operators, see
Wigner (\citeyear{wigner51}, \citeyear{wigner52}). Later, in \citet{Wigner57},
he turned to models of random matrices that are by now
standard, such as the Gaussian orthogonal ensemble (GOE).
In this model, one fills an $n\times n$ matrix $M$ with
independent standard normal random variables, then
symmetrizes it to get
\[
A=\frac{M+M^{T}}{\sqrt2}.
\]
The~Wigner semicircle law is the limit of the empirical
distribution of the eigenvalues of the matrix $A$. However,
Wigner's main interest was the local behavior of the eigenvalues,
in particular the repulsion between them.
He examined the asymptotic probability of having no eigenvalue in a
fixed interval of size $\lambda$ for $n\to\infty$ while the spectrum
is rescaled to have an average eigenvalue spacing $2\pi$.
%
Wigner's prediction for this probability was
\[
p_\lambda= \exp\bigl(- \bigl(c +o(1)\bigr)\lambda^2\bigr),
\]
where this is a $\lambda\to\infty$ behavior.
%
This
rate of decay is in sharp contrast with the exponential
tail for gaps between Poisson points; it is one
manifestation of the more organized nature of the random
eigenvalues. Wigner's estimate of the constant $c$,
$1/(16\pi)$, later turned out to be inaccurate. \citet{Dy62}
improved this estimate to
%
\begin{equation}\label{mainform}
p_\lambda
=\bigl(\kappa_\beta+o(1)\bigr)\lambda^{\gamma_\beta} \exp\biggl(-
\frac{\beta}{64}\lambda^2+\biggl(\frac{\beta}{8}-\frac14
\biggr)\lambda\biggr),
\end{equation}
where $\beta$ is a new parameter introduced by noting that the
joint eigenvalue density of the GOE is the $\beta=1$ case of
%
\begin{equation}
\label{betadens}
\frac{1}{Z_{n,
\beta}}
e^{ - \beta\sum_{k=1}^n \lambda_k^2/4}   \prod_{j < k} |
\lambda_j - \lambda_k |^{\beta}.
\end{equation}
The~family of distributions defined by the density \eqref{betadens} is
called the $\beta$-ensemble.
Dyson's computation of the exponent $\gamma_\beta$, namely
$\frac{1}{4}(\frac{\beta}{2}+\frac{2}{\beta}+ 6)$, was
shown\vspace*{1pt} to be
slightly incorrect. Indeed, \citet{dCM73} gave more substantiated
predictions that $\gamma_\beta$ is equal to $-1/8, -1/4$ and $-1/8$
for values $\beta=1, 2$ and~4, respectively. 
Mathematically precise proofs for the $\beta=1,2$ and 4 cases were
later given by several authors: \citet{Wi96} and \citet{DIZ96}. Moreover,
the value of $\kappa_\beta$ and higher-order asymptotics were also
established for these specific cases by \citet{Kr04},
\citet{Ehr06} and
\citet{DIKZ07}.
The~problem of determining the asymptotic probability
of a large gap naturally arises in other random matrix models as well.
In the physics literature, \citet{ChenMan} treat the problem of the
$\beta$-Laguerre ensemble at the edge.

Our main theorem gives a mathematically rigorous version of Dyson's
prediction for general $\beta$ with a corrected exponent $\gamma
_\beta$.

\begin{theorem}\label{mainthm}
The~formula \eqref{mainform} holds
with a positive $\kappa_\beta$ and
\[
\gamma_\beta=\frac{1}{4}\biggl(\frac{\beta}{2}+\frac{2}{\beta}-
3\biggr).
\]
\end{theorem}

The~proof is based on the Brownian carousel, a geometric
representation of the $n\to\infty$ limit of the eigenvalue
process. We first introduce the \textit{hyperbolic carousel}. Let:
\begin{itemize}
\item$b$ be a path in the hyperbolic plane,
\item$z$ be a point on the boundary of the hyperbolic plane
and
\item$f\dvtx \mathbb R_+\to\mathbb R_+$ be an integrable function.
\end{itemize}
To these three objects, the hyperbolic carousel associates a
multi-set of points on the real line defined via its counting
function $N(\lambda)$ taking values in $\mathbb
Z\cup\{-\infty,\infty\}$. As time increases from $0$ to $\infty$,
the boundary point $z$ is rotated about the center $b(t)$ at
angular speed $\lambda f(t)$. $N(\lambda)$ is defined as the
integer-valued total winding number of the point about the moving
center of rotation.

\begin{figure}

\includegraphics{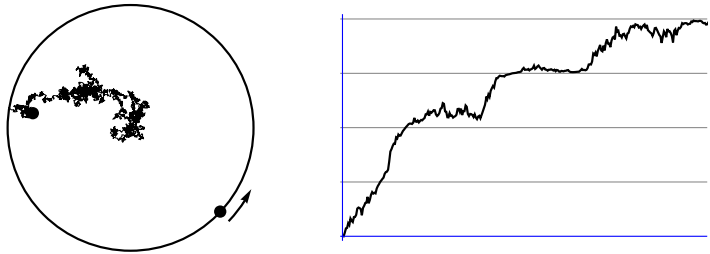}

\caption{The~Brownian carousel and the winding angle $\alpha_\lambda$.}\label{fig1}
\end{figure}

The~\textit{Brownian carousel} is defined as the
hyperbolic carousel driven by hyperbolic Brownian motion
$b$ (see Figure \ref{fig1}). It is connected to random matrices via
the following theorem.%

\begin{theorem}[{[\citet{BVBV}]}]
Let $\Lambda_n$ denote the point
process given by \eqref{betadens}, and let $\mu_n$ be a sequence
so that  $n^{{1}/{6}}(2\sqrt{n}-|\mu_n|)\to\infty$. Then we
have the following convergence in distribution:
%
\begin{equation}\label{e_limit_eredeti}
\sqrt{4n-\mu_n^2}(\Lambda_n-\mu_n)\quad  \Rightarrow\quad \mathrm{Sine}_\beta,
\end{equation}
where $ \mathrm{Sine}_\beta$ is the discrete point process given by
the Brownian carousel with parameters
%
\begin{equation}
f(t)=\frac{\beta}{4}e^{-\beta
t/4} \label{e:drift}
\end{equation}
and arbitrary $z$.
\end{theorem}

\begin{remark}
The~semicircle law shows that most points in $\Lambda_n$ are in the
interval $[-2\sqrt{n},2\sqrt{n}]$.
The~discrete point process $\Lambda_n$ has two kind of point process
limits, one near the edges of this interval and another in the bulk.
The~condition on the parameter $\mu_n$
means that we get a bulk-type scaling limit of $\Lambda_n$. The~scaling factor in (\ref{e_limit_eredeti}) is the natural choice in
view of the Wigner semicircle law in order to get a point process with
average density $1/(2\pi)$.
The~limiting point process for the edge-scaling case have been obtained
by \RRVlong.
\end{remark}

The~Brownian carousel description gives a simple way to analyze
the limiting point process. The~hyperbolic angle of the rotating
boundary point as measured from $b(t)$ follows the following
coupled one-parameter family of stochastic differential equations
%
\begin{equation}\label{e_ssefenetudja}
d\alpha_\lambda= \lambda f\, dt + \operatorname{Re}
\bigl((e^{-i\alpha_\lambda}-1)\,dZ\bigr),\qquad   \alpha_\lambda(0)=0,
\end{equation}
driven by a two-dimensional standard Brownian motion and $f$ given in
(\ref{e:drift}). For a single
$\lambda$, this reduces to the one-dimensional stochastic
differential equation
%
\begin{equation}\label{eq:sse}
d\alpha_\lambda= \lambda f\, dt + 2
\sin(\alpha_\lambda/2)\,dW, \qquad \alpha_\lambda(0)=0,
\end{equation}
which converges as $t\to\infty$ to an integer multiple
$\alpha_\lambda(\infty)$ of $2\pi$. In particular, the
number of points of the point process $\mathrm{Sine}_\beta$ in
$[0,\lambda]$ has the same distribution as
$\alpha_\lambda(\infty)/(2\pi)$ and $p_\lambda$ is equal to the
probability that $\alpha$ converges to $0$ as $t\to\infty$.
See \citet{BVBV} for further details.

In the analysis of equation \eqref{eq:sse}, it helps to
remove the space dependence from the diffusion coefficient
by a change of variables $X(t)=\log(\tan(\alpha(t)/4))$. The~diffusion $X$ satisfies the stochastic differential
equation
%
\begin{eqnarray}\label{eq:logtan}
dX= \frac{\lambda}{2}f\cosh X \,dt+\frac{1}{2}\tanh X
\,dt+dB,\qquad   X(0)=-\infty.
\end{eqnarray}
In \citet{BVBV}, equations \eqref{eq:sse} and
\eqref{eq:logtan} were used to identify the leading term
in the asymptotic expansion of $p_\lambda$ in
\eqref{mainform}. The~proof of Theorem \ref{mainthm}
requires a more careful analysis of equation
\eqref{eq:logtan}.

In Lemma \ref{l_unique}, we will show that for any initial
condition $X(0)=x\in[-\infty,\infty)$ there is a unique
solution of the equation given in (\ref{eq:logtan}) and
the desired gap probability $p_\lambda$ may be written in terms of a
passage probability for this process. Namely,
$p_\lambda=p_\lambda(-\infty)$ where
%
\begin{eqnarray}\label{eq:gapX}
p_\lambda(x)&:=&\mathbf{P}\bigl(
\mbox{$X(t)$ is finite for all $t>0$ and}\nonumber\\[-8pt]\\[-8pt]
&&\quad\mbox{it does not converge to $+\infty$ as
$t\to\infty$}\bigr).\nonumber
\end{eqnarray}
A time
shift of equation \eqref{eq:logtan} only changes the
parameter $\lambda$ and the initial condition. This, together
with the Markov property of the diffusion $X$, shows that
with \mbox{$T=\frac4{\beta}\log\lambda$} we have
\begin{eqnarray*}
p_\lambda&=&\mathbf{E}[ \mathbf{1}\{ \mbox{$X(t)$ is
finite for all
$0<t\le T$}\} \cdot p_1(X(T))].
\label{eq:condprob}
\end{eqnarray*}
Our main tool is the Cameron--Martin--Girsanov formula, which
allows one to compare the measure on paths given by two
diffusions. If we knew the conditional distribution of the
diffusion $X$ under the event that it does not blow up,
then we could use the Cameron--Martin--Girsanov formula to compute
$p_\lambda$
explicitly. While we cannot do this, the next best option
is to find a new diffusion $Y$ which approximates this
conditional distribution. The~density (i.e., the
Radon--Nikodym derivative) of the path measures given by $Y$
with respect to the measure given by $X$ will be close to
the right-hand side of $\eqref{mainform}$. Our strategy for
finding $Y$ is described in Section \ref{sec:Y}.
Cameron--Martin--Girsanov techniques have been used to obtain tail
asymptotics for the tail asymptotics of the ground state of the random
Hill's equation, see \citet{CRR06}.

In Section \ref{sec:proof}, we will present a coupling of
the transformed processes that enables us to show that the
asymptotics is precise up to and including the constant
term $\kappa_\beta$. The~term $\kappa_\beta$ is then
identified as the expectation of a functional of a certain
limiting diffusion.

\subsection*{Open problem}
Give an explicit expression for
$\kappa_\beta$ for general values of $\beta$.

The~known values of $\kappa_\beta$ are
\[
\kappa_1=2^{13/24} e^{({3}/{2}) \zeta'(-1)},\qquad   \kappa
_2=2^{7/12} e^{3 \zeta'(-1)},\qquad   \kappa_4=2^{-13/12} e^{({3}/{2}) \zeta'(-1)},
\]
where $\zeta'(-1)$ is the coefficient of the linear term in the
Laurent series of the Riemann-$\zeta$ function at $-1$.

A natural generalization of Theorem \ref{mainthm} would be to consider
the asymptotic probability that there are exactly $k$ eigenvalues in a
large interval $[0,\lambda]$. This probability is usually denoted by
$E_\beta(k;\lambda)$ in the literature. For $\beta=1,2$ and 4 the
following large $\lambda$ asymptotics was obtained by \citet{BTW}:
\[
\log E_\beta(k;\lambda)=\log E_\beta(0;\lambda)+\frac{k \beta}{4}
\lambda+\frac{k}{2}\biggl(1-\frac{\beta}{2}-\frac{k \beta
}{2}\biggr)
\log\lambda+c_\beta+o(1)
\]
with explicit constants $c_\beta$. [See also \citet{TrWid}.] We
believe that our methods can be used to extend the previous asymptotics
for general values of $\beta$.

The~rest of the paper is organized as follows. In the next
section, we justify the $p_\lambda=p_\lambda(-\infty)$ and
give some preliminary estimates on the probability $p_1(r)$
appearing in \eqref{eq:gapX}. Section \ref{sec:CMG}
presents the version of the Cameron--Martin--Girsanov formula
that we need. Section~\ref{sec:Y} describes the strategy
for finding $Y$ and Section~\ref{sec:proof} builds on these
sections to complete the proof of the main theorem.

\section{Preliminary results}

First, we formally verify the connection between the gap
probability $p_\lambda$ and the diffusion given in
\eqref{eq:logtan}.

\begin{lemma}\label{l_unique}
The~diffusion \eqref{eq:logtan} has a unique solution for any initial
condition $X(0)=x\in[-\infty,\infty)$ and $p_\lambda=p_\lambda
(-\infty)$.
\end{lemma}

\begin{pf}
The~change of variables function $\log\tan(\cdot/4)$ is
one-to-one on $(0,2\pi)\to\mathbb{R}$. Therefore, even with
$-\infty$ initial condition, the diffusion $X$ is well
defined and has a unique solution until $\alpha$ reaches
$2\pi$, when it blows up. We define $X(t)=\infty$ after
this blowup.

Note that for $\lambda>0$ the solution of equation
\eqref{eq:sse} is always monotone increasing at multiples
of $2 \pi$. See Section 2.2 in \citet{BVBV} for more details.
So, if $\alpha(t) \to0$ as $t\to\infty$ then
$0<\alpha(t)<2\pi$ for all $t>0$. This means that $X(t)$ is
finite for all $t>0$ and $X(t)$ cannot converge to
$\infty$ which proves $p_\lambda=p_\lambda(-\infty)$.
\end{pf}

Next, we prove a preliminary estimate on the blowup
probability of the diffusion~\eqref{eq:logtan}.

\begin{lemma}\label{l:hapci}
Recall that $p_1(x)$ is the probability that the diffusion~\eqref{eq:logtan} with \mbox{$\lambda=1$}
and initial condition
$X(0)=x$ does not blow up in finite time and does not
converge to $+\infty$ as $t \to\infty$. We have
\[
0<p_1(x)\le c_\beta\exp\biggl( -\frac{\beta}{60}
 e^{x}\biggr).
\]
\end{lemma}

\begin{pf}
%
For the upper bound, we first assume that $x>4$. Consider the diffusion
%
\begin{equation}\label{e_YY}
d R = \frac{\beta}{16}  e^{R-\beta t/4} \, dt + dB, \qquad   R(0)=x.
\end{equation}
This has the same noise term as $X$. The~drift term of $R$ is $\frac
{\beta}{16}  e^{x-\beta t/4 }$, which
is\vspace*{-1pt} dominated by the drift term $\frac{f(t)}{2}\cosh(x)$ of $X$ when
$x$ is nonnegative.
Thus, while $R$ stays positive,
we have $R\le X$. This means that for every $t>0$ we have
\begin{eqnarray}\label{Rprob}
\nonumber
p_1(x) &\le& \mathbf{P}( \mbox{$X$
does not blow up before time } t )
\\[-8pt]\\[-8pt]
&\le&\mathbf{P}\Bigl(\min_{s\in[0,t]}R(s)<0 \mbox{ or $R$ does not
blow up before time } t\Bigr) .\nonumber
\end{eqnarray}
The~difference $Z= R-B$ satisfies the ODE
%
\begin{equation}\label{e_ZZ}
e^{-Z}\,dZ= \frac{\beta}{16}  e^{B-\beta/4   t} \, dt , \qquad   Z(0)=x.
\end{equation}
Integration gives
\[
e^{-x}-e^{-Z(t)}= \frac{\beta}{16} \int_0^t e^{B(s)-\beta/4
  s} \, ds.
\]
This shows that $Z$ is increasing in $t$, in particular
$Z(t)\ge x$. So if $\min_{[0,t]} R<0$ then
\[
\min_{[0,t]} B <-x.
\]
Furthermore, if
%
\begin{equation}\label{blowup}
e^{-x}<\frac{\beta}{16} \int_0^t e^{B(s)-\beta/4   s}
 \,ds
\end{equation}
then $R$ blows up before time $t$. This certainly happens if the
minimum of $B$ on the interval $[0,t]$ is not sufficiently small.
More precisely, if
%
\begin{equation}\label{kuka2} \frac{e^{-b} }{4} (1-e^{\beta t/4})>e^{-x}
\end{equation}
and $\min_{[0,t]} B>-b$ then \eqref{blowup} follows. So if
$b<x$, and \eqref{kuka2} holds, then the right-hand side of
\eqref{Rprob} can be bounded above by
\[
P\Bigl(\min_{[0,t]} B<-b\Bigr)=P\bigl(|B(t)|>b\bigr)\le
\frac{\sqrt{t}}{b}  e^{-{b^2}/{(2t)}}.
\]
We set
\[
t= \frac{16}{\beta}  e^{2-x}, \qquad   b=
\frac{4e}{\sqrt{30}}<2.
\]
As $x>4$, both $b<x$ and \eqref{kuka2} are satisfied and we get the
upper bound
\[
p_1(x)\le\frac{\sqrt{t}}{b}  e^{-{b^2}/{(2t)}}<c_\beta
e^{-({\beta}/{60}) e^x}
\]
with $c_\beta=\sqrt{{30/\beta}}$. The~upper bound for all values of
$x$ now follows by changing the constant $c_\beta$ appropriately.

For the lower bound note that since the $\operatorname{Sine}_{\beta}$
process is
discrete and translation invariant in distribution, there exists
$\nu\in(0,1)$ so that $p_\nu=p_{\nu}(-\infty)>0$. By the Markov
property, we
have
\[
p_\nu= \int_{-\infty}^{\infty} K_{0,1}(-\infty,dx)p_{\nu
e^{\beta/4}}(x),
\]
where $K_{s,t}(y,dx)$ is the transition kernel of the Markov process
$X$ with parameter $\lambda=\nu$. This implies that for
some $x_0\in\mathbb{R}$ we have
\[
p_{\nu e^{\beta/4}}(x_0)>0.
\]
Consider the process $X$ started at $x$ with parameter
$\lambda=1$. The~Markov property applied at time
$t_0=1-\frac{4}{\beta}\log\nu$ and the monotonicity of
$p_\lambda(x)$ in $x$ implies
\[
p_1(x) \ge P\bigl( X(t_0)<x_0\bigr) p_{\nu
e^{\beta/4}}(x_0)
>0,
\]
since $P(X(t)<x)$ is positive for all $x\in\mathbb{R}$ and $t>0$.
%
%
\end{pf}

\section{The~Cameron--Martin--Girsanov formula} \label{sec:CMG}

Our main tool will be the following version of the
Cameron--Martin--Girsanov formula. Here, we allow diffusions to blow
up to $+\infty$ in finite time, in which case they are required to
stay there forever after.
%

\begin{proposition}\label{pgirsanov}
Consider the following stochastic differential equations:
%
\begin{eqnarray}
dX&=&g(t,X)\, dt+dB,\qquad  \lim_{t\to0} X(t)=-\infty,\label{sde1}\\
dY&=&h(t,Y) \,dt+d\tilde B,\qquad  \lim_{t\to0} Y(t)=-\infty
\label{sde2}
\end{eqnarray}
on the interval $(0,T]$ where $B, \tilde B$ are standard Brownian
motions. Assume that \eqref{sde1} has a unique
solution $X$ in law taking values in $(-\infty, \infty]$.

Let
%
\begin{equation}\label{e:Girsx}
\qquad G_s=G_s(X)=\int_0^s\bigl( h(t,X)-g(t,X)\bigr)\, dX-\frac12 \int_0^s
\bigl(h(t,X)^2-g(t,X)^2\bigr)\, dt
\end{equation}
and assume that:
\begin{longlist}
\item[\textup{(A)}]\label{c1} $ g^2-h^2 $ and $g-h$ are bounded when $x$ is bounded above.
(Then $G_s$ is almost surely well defined when $X_s$ is finite.)

\item[\textup{(B)}]\label{c2} $G_s$ is bounded above by a deterministic constant.

\item[\textup{(C)}]\label{c3} $G_s\to-\infty$ when $s\uparrow\tau$ if $X$ hits
$+\infty$ at time $\tau$. In this case, we define $G_s:=-\infty$
for $s\ge\tau$.
\end{longlist}

Consider the
process $\tilde Y$ whose density with respect to the distribution
of the process $X$ is given by $e^{G_T}$. Then $\tilde Y$
satisfies the second SDE \eqref{sde2} and never blows up to
$+\infty$ almost surely. Moreover, for any nonnegative function
$\varphi$ of the path of $X$ that vanishes when $X$ blows up we
have
%
\begin{equation}\label{xtoy}
\mathbf{E} \varphi(X)=\mathbf{E}\bigl[\varphi(Y)
e^{-G_T(Y)}\bigr].
\end{equation}
\end{proposition}

\begin{remark}
There exist several versions of the Cameron--Martin--Girsanov
formula for exploding diffusions [e.g.,~\citet{McK}, Section 3.6]. As
we did not find one in the literature which could be directly applied
to our case, we sketch the proof below.
\end{remark}

\begin{pf*}{Proof of Proposition \ref{pgirsanov}}
We follow the standard proof of the Girsanov
theorem.

First, we show that $G_s$ is well defined for finite $X_s$. From
condition \textup{(A)}, it follows that the second integral is
well defined. The~first integral can be written as
\[
\int_0^s (h-g)\,dB+\int_0^s (h-g)g\, dt
\]
which is well defined since $(h-g)$ and
$2(h-g)g=(h^2-g^2)-(h-g)^2$ when their argument $x$ is bounded
above.

Next, we show that $M_s=e^{G_s}$ is a bounded martingale.
This is clear after the hitting time $\tau$ of $X$ of
$+\infty$, if such time exists. Before this time, $G_s$ is
a semimartingale, and so is $M_s$. It\^{o}'s formula gives
\[
dM= (h-g)M\, dB
\]
so that the drift term of $M$ vanishes. So $M$ is a local
martingale which is bounded, so it has to be a martingale.

The~rest of the proof is standard and is outlined as follows. Set
\[
\tilde B_s=X_s-\int_0^s h \, dt=B_s-\int_0^s (h-g)\,dt.
\]
It suffices to show that $\tilde B$ is a Brownian motion with
respect to the new measure with density $M_T$. This follows from
L\'{e}vy's criterion [\citet{KarShr}, Theorem 3.3.16] if $\tilde B$
and $\tilde B^2-s$ are local martingales. Since $M$ is a
martingale, it suffices to show that $\tilde BM$ and
$(\tilde B^2-s)M$ are local martingales with respect to the old
measure, which is just a simple application of It\^{o}'s formula.

The~identity \eqref{xtoy} is just a version of the change of
density formula.
\end{pf*}

\section{Construction of the diffusion $Y$} \label{sec:Y}

In this section, we will create
a diffusion which approximates the
conditional distribution of the
diffusion $X$ under the event that it does not blow up.
We will construct a drift function $h(t,x)$ for
which the diffusion $Y$
%
\begin{equation}\label{eq:Y}
dY=h(t,Y)\, dt+dB_t,\qquad   Y(0)=-\infty,
\end{equation}
is well defined, a.s.~finite for $t>0$ and the (formal)
Radon--Nikodym derivative $e^{G_T}$ with $G_T$ defined in (\ref{e:Girsx}) is almost equal to the
right-hand side of equation \eqref{mainform} with the appropriate
$\gamma_\beta$.

\begin{lemma}\label{l_Y}
For the diffusion \eqref{eq:logtan}, $\lambda>1$ and
$T=\frac4\beta\log\lambda$ there exists a function
$h(t,x)$ so that conditions \textup{(A)}--\textup{(C)} of
Proposition \textup{\ref{pgirsanov}} hold, and $G_T$ has the
following form:
%
\begin{eqnarray}\label{eq:Girsan}
-G_T(X)&=& -\frac{\beta}{64}  \lambda^2
+\biggl(\frac{\beta}{8}-\frac14
\biggr) \lambda+\frac18\biggl(\beta+\frac4{\beta}-6 \biggr)\log
\lambda\nonumber
\\
&&{}+\frac{ \beta}{8} e^{X(T)}+\biggl(2-\frac{\beta
}{2}\biggr) X(T)^+
+\omega(X(T))\\
&&{}+\int_{0}^T \phi\bigl(T-t,X(t)\bigr) \,dt .\nonumber
\end{eqnarray}

Here, the function $\omega$ is bounded and continuous, $\phi$
is continuous and $|\phi(t,x)|\le\tilde\phi(t)$ with
$\int^{\infty}_0 \tilde\phi(t)\, dt<\infty$. The~functions $\omega$
and $\phi,\tilde\phi$ may depend on the parameter~$\beta$, but not
on $\lambda$.

The~function $h$ will have the following form:
%
\begin{equation}\label{eq:www}
h(t,x)=-\frac{\lambda}2 f \sinh(x)+ h_0(t,x),
\end{equation}
where $| h_0(t,x)|<c$ if $0\le t\le T$. The~constant $c$ depends only
on $\beta$.
\end{lemma}

\begin{pf}{\textit{Construction of the function} $h$}.
 Given an explicit
formula for $h$ it would not be hard to check that ${G_T}$ has the
desired form. However, we would like to present a way one
can find the appropriate drift function. This will provide
a better understanding of the form of the resulting $h$.

We will use the definition
\[
-G_s(X)=\int_0^s \bigl(g(t,X)-h(t,X)\bigr)\, dX+\frac12 \int_0^s
\bigl(h^2(t,X)-g^2(t,X)\bigr)\,dt,
\]
where
\[
g=g_1+g_2, \qquad  g_1(t,x)=\frac{\lambda}{2}f(t)\cosh x,\qquad
g_2(t,x)=\frac{1}{2}\tanh x.
\]
Our goal is to find the appropriate drift term $h$
in a way that the diffusion $Y$ will approximate the
conditional distribution of $X$ given that it does not blow
up in the interval $[0,T]$. We will do this term by term,
starting with the highest order; toward this end we write
$h=h_1+h_2+h_3+h_4$. We set
%
\begin{equation}\label{eq:h1}
h_1(t,x)= -\frac{\lambda}{2} f(t)\sinh(x)
\end{equation}
as this yields the nice cancelation
\[
h_1^2-g_1^2= \frac{\lambda^2}{4} f(t)^2
\sinh^2(X)-\frac{\lambda^2}{4} f(t)^2
\cosh^2(X)=-\frac{\lambda^2}{4} f(t)^2
\]
in the main terms of $h^2-g^2$. In addition, if the
remaining term $h_2+h_3+h_4$ is bounded, then it will be
easy to show that conditions \textup{(A)}--\textup{(C)} of
Proposition \ref{pgirsanov} are satisfied. This will be done at the
end of the proof.

The~contribution of the drift terms $h_1$ and $g_1$ to the
stochastic integral part of~$-G_s$ is given by
%
\begin{equation}\label{eq:dXterm}
\frac{\lambda}{2} \int_0^s f(t)\bigl(\cosh(X)-\sinh(X)\bigr)\, dX=
\frac{\lambda}{2} \int_0^s e^{X}f(t)\,dX.
\end{equation}
Our main tool for evaluating integrals with respect to $dX$ is the
following version of It\^{o}'s formula. Let $a,b$ be
continuously differentiable functions and let $\tilde a$
denote the antiderivative of $a$. Then
%
\begin{equation}\label{eq:ito}
a(t)b(X) \,dX = d(a(t)\tilde b(X)) - a'(t)\tilde
b(X) \,dt-\tfrac{1}{2} a(t)b'(X)\, dt.
\end{equation}
%
Since $f'(t)=-\beta/4 f(t)$, and $X(0)=-\infty$, this
formula gives
%
\begin{eqnarray}
\frac{\lambda}{2} \int_0^s f(t)e^{X}
\,dX = \frac{\lambda}{2}
f(s)e^{X(s)}+\frac{\lambda}{2}\biggl(\frac{\beta}{4}-\frac12
\biggr) \int_0^s e^{X} f\, dt.\label{eq:tag1}
\end{eqnarray}
Next, we would like to choose $h_2$ in \eqref{eq:h1} so that
the integral term in the right-hand side of \eqref{eq:tag1}
simplifies. More precisely, since we expect the diffusion
$X$ to be near~0 most of the time, we would like to replace
the term $e^X$ by $1$. The~plan is to use the cross term
$\int h_1h_2 \, dt$ in the $\frac{1}{2}\int h^2 \, dt$ term
of $G$ to do this. Namely, we would like to have
%
\begin{equation}\label{eq:tag2}
 h_1h_2=\frac{\lambda}{2}
\biggl(\frac
{\beta}{4}-\frac12 \biggr)
(1-e^x) f.
\end{equation}
The~solution for (\ref{eq:tag2}) is given by
%
\begin{equation}
h_2(t,x)=\biggl(\frac{\beta}{4}-\frac12 \biggr)
\bigl(1+\tanh(x/2)\bigr).\label{h2}
\end{equation}
We will choose the next term, $h_3$, so that the cross term
$\int h_1h_3\, dt$ in $\frac{1}{2}\int h^2\, dt$ cancels the
cross term $-\int g_1g_2$ in $-\frac12 \int g^2  \, dt$.
This leads to the equation
\[
h_1h_3=g_1g_2=\frac{\lambda}{2} f(t)\cosh(x)\cdot
\frac12\tanh(x),
\]
which gives
%
\begin{equation}
h_3(t,x)=-\tfrac12.\label{h3}
\end{equation}
Collecting all our
previous computations, we get
%
\begin{eqnarray}
-G_s&=&\frac{\lambda}{2} f(s)e^{X(s)}
-\frac{\lambda^2}{8}\int_0^s f^2 \, dt+\lambda\biggl(\frac{\beta
}{8}-\frac
14 \biggr) \int_0^s f  \,dt\nonumber\\
&&{}+\frac12 \int_0^s 2h_1h_4+(h_2+h_3+h_4)^2 -g_2^2 \,dt\label
{eq:Girs}\\
&&{}-\int_0^s h_4\, dX+\int_0^s(g_2-h_2-h_3)\,  dX.\nonumber
\end{eqnarray}
%
The~integrand $u=g_2-h_2-h_3$ in the last integral of
\eqref{eq:Girs} has antiderivative
%
\begin{equation}\label{eq:qqq}
\tilde u(x)= \biggl(1-\frac\beta4 \biggr)x + \biggl(1-
\frac{\beta}2\biggr) \log\cosh(x/2) +
\frac12 \log\cosh x.
\end{equation}
By It\^{o}'s formula, $\int_0^s u(X)\, dX-\tilde u(X)
|_0^s$ is given by
\begin{eqnarray}\label{upi}\nonumber
-\frac12\int_0^s u'(X) \,dt&=&-\frac12
\int_0^s \biggl[ \frac{2-\beta}{8} \operatorname{sech}(X/2)^2+\frac
12 \operatorname{sech}
(X)^2\biggr]\, dt
\\[-8pt]\\[-8pt]
&=& \frac{\beta-6}{16} s+ \int_0^s
\biggl[ \frac{( 2 - \beta)}{16}
 {\tanh(X/2)}^2 +\label{eq:tag3} \frac14 {{\tanh(X)}^2}
\biggr]\, dt.\nonumber
\end{eqnarray}
Note that
\[
\lim_{x\to-\infty} \tilde u(x)=\frac{\beta-3}{2}
\log2=c_1.
\]
Substituting this computation for the last integral and
expanding $(h_2+h_3+h_4)^2$, we can rewrite \eqref{eq:Girs}
as follows:
%
\begin{eqnarray}\label{eq:Girs1}\nonumber
-G_s&=&-\frac{\lambda^2}{8}\int_0^s f^2 \, dt+\lambda
\biggl(\frac{\beta}{8}-\frac14 \biggr) \int_0^s f \,dt\\
&&{}+\biggl(
\frac12\biggl(\frac{\beta}{4}-1 \biggr)^2+\frac{\beta-6}{16}
\biggr)s\nonumber
\\
&&{}+\frac{\lambda}{2} f(s)e^{X(s)}+\tilde u(X(s))-c_1 \\
&&{}+ \frac{1}{2}\int_0^s\bigl( 2h_1 h_4
+2(h_2+h_3)h_4+h_4^2\bigr)\,
 dt\nonumber\\
 &&{}-\int_0^s h_4\, dX
+\int_0^s \eta(X(t))\, dt \nonumber.
\end{eqnarray}
The~coefficient of $s$ in the first line of
\eqref{eq:Girs1} comes from the first term on the right-hand side in
\eqref{upi} and the constant term of $(h_2+h_3)^2/2$. The~function $\eta$ collects the terms from the integrand in
\eqref{upi}, the terms $(h_2+h_3)^2/2$ with the constant
term $(\beta/4-1)^2/2$ removed, and $-g_2^2/2$. More
explicitly, we have
\begin{eqnarray}
\eta(x)=\frac{( 8
- 6 \beta+ {\beta}^2 )}{32}
\bigl( 2 \tanh(x/2) + {\tanh(x/2)}^2 \bigr) +
\frac{1}{8}(\tanh x)^2.
\nonumber
\end{eqnarray}
The~function $\eta(x)$ contributes to an error term that
needs to be controlled, but whose precise value does not
influence our final result. Now, we are ready to set the
value for $h_4$: we will choose it in a way that the cross
term $\int h_1 h_4\, dt$ in \eqref{eq:Girs1} will cancel the
integral $\int\eta \, dt$. This gives $h_4=-\eta/h_1$,
that is,
%
\begin{equation}
h_4= \frac{2}{\lambda f(t)} \frac{\eta(x)}{\sinh(x)}.
\label{eq:gam2}
\end{equation}
The~function $h_4$ is a product of a function of $t$ and a
function of $x$. It\^{o}'s formula~\eqref{eq:ito}, with the
notation $\tilde h_4(t,x)=\int_0^x h_4(t,y)\,dy$ yields the
evaluation of the stochastic integral in \eqref{eq:Girs1}:
\begin{eqnarray*}
-\int_0^s h_4 \,dX = -\tilde
h_4(s,X(s))+\frac{\beta}{4}\int_0^s \tilde h_4\, dt+\frac12
\int_0^s \partial_x h_4\, dt.
\end{eqnarray*}
Plugging this into \eqref{eq:Girs1} and simplifying the
deterministic terms in the first line of~\eqref{eq:Girs1},
we arrive at
%
\begin{eqnarray}\label{eq:Girs2}
-G_s&=&
-\frac{\lambda^2}{8}\int_0^s f^2 \, dt+\lambda\biggl(\frac{\beta
}{8}-\frac
14 \biggr) \int_0^s f \, dt+\frac1{32}(\beta^2+12{\beta}+8
) s\nonumber
\\
&&{}+\lambda e^{-\beta/4 s} \frac{ \beta}{8} e^{X(s)}+\tilde
u(X(s))-c_1-\tilde h_4(s,X(s))\\
&&{}+ \int_0^s\biggl(2(h_2+h_3)h_4+h_4^2+
\frac{\beta}{4}\tilde h_4+\frac12 \partial_x h_4 \biggr)\,
 dt\nonumber.
\end{eqnarray}
%
Note that $h_2$ and $h_3$ do not depend on $t$ and are bounded by an
absolute constant. The~functions $h_4, \tilde
h_4,
\partial_x h_4$ are all bounded by a constant times $1/(\lambda
f(t))=\frac{16}{\beta^2}
f(T-t)$, which itself is bounded by a constant not
depending on $\lambda$ as long as\vspace*{1pt} $0\le t\le T$. Thus, we
can rewrite the integrand in \eqref{eq:Girs2} as
%
\begin{equation}\label{e:tag77}
\int_{0}^s \phi\bigl(T-t,X(t)\bigr)\, dt
\end{equation}
with a continuous function $\phi$ which does not depend
on $\lambda$ and satisfies $|\phi(t,x)|\le\tilde\phi(t)$ with
$\int^{\infty}_0 \tilde\phi(t) \,dt<\infty$. Using \eqref{eq:qqq}
and the fact that
$\log\cosh x-|x|$ is bounded, the terms in the second line
of \eqref{eq:Girs2} can be written
as
%
\begin{equation}
\label{e:tag88} \biggl(2-\frac\beta2 \biggr)X(s)^+ +\lambda
e^{-\beta/4 s} \frac{\beta}{8} e^{X(s)}
+\omega_0(X(s))-\tilde h_4(T,X(T))
\end{equation}
with a bounded and continuous $\omega_0$. This concludes
the construction of the function~$h$. In order to get the
expression (\ref{eq:Girsan}) for $-G_T$, we first plug in
$s=T$ into~\eqref{eq:Girs2}. Then the first line gives
\[
-\frac{\lambda^2 \beta}{64}(1-\lambda^{-2})+\lambda
\biggl(\frac{\beta}{8}-\frac14 \biggr)
(1-\lambda^{-1})+\frac14\biggl(\frac\beta2+\frac
2{\beta}-3 \biggr)\log\lambda,
\]
and by (\ref{e:tag88}) the second line transforms to
\[
\biggl(2-\frac\beta2 \biggr)X(s)^+ + \frac{\beta}{8}
e^{X(s)} +\omega_0(X(T))-\tilde h_4(T,X(T)).
\]
Note that the expression $\tilde h_4(T,x)$ does not depend
on $T$ and is bounded. This proves that $-G_T$ is in the
desired from (\ref{eq:Girsan}).

Now, we are ready to check that the
proposed choice of $h$ satisfies all the needed conditions
(A)--(C).

\textit{Condition} \textup{(A)}.
%
As $x\to
-\infty$, we have
\[
g(t,x)=\tfrac14 \lambda f e^{-x}-\tfrac12+\hat g(t,x),\qquad
h(t,x)=\tfrac14 \lambda f e^{-x}-\tfrac12+\hat h(t,x),
\]
where $|\hat g|<c  e^{x}$ and $|\hat h|<c' e^{x}$ with
constants that only depend on $\beta$ if $0\le t \le T$.
From this, it follows that
$g-h$ and $g^2-h^2$ are both bounded if $x$ is bounded from above.

\textit{Condition} \textup{(B)}. We need that \eqref{eq:Girs2}
is bounded from below if $0\le s\le T$. The~integrals in
the first line are bounded by a constant depending
on $\lambda$ and $\beta$ only. The~same is true for the integral
in the last line, see \eqref{e:tag77} and the discussion
around it. Thus, we only need to deal with the evaluation
terms of the second line. By \eqref{e:tag88}, we just need
to show that
%
\begin{equation}\label{verify}
\biggl(2-\frac\beta2 \biggr)X(s)^+ +\lambda e^{-\beta/4 s}
\frac{\beta}{8} e^{X(s)}
\end{equation}
is bounded from below. Since $ s \le T=\frac{4}{\beta} \log{\lambda}$,
we get that (\ref{verify}) is bounded from below by
\[
\biggl(2-\frac\beta2
\biggr) X(s)^+ + \frac{\beta}{8} e^{X(s)}
\]
which in turn is bounded from below by a constant depending
only on $\beta$.

\textit{Condition} \textup{(C)}. This follows the same way: one
only needs to check the behavior of \eqref{verify} as
$s$ converges to the hitting time of $\infty$. This expression
converges to $\infty$ as $X(s)\to\infty$ which means that
$G_s \to-\infty$.
\end{pf}

\section{The~proof of the main theorem} \label{sec:proof}

We are ready to prove Theorem~\ref{mainthm}.

\begin{pf*}{Proof of Theorem \ref{mainthm}}
Lemma \ref{l_unique} gives $p_\lambda=p_\lambda(-\infty)$,
where
\begin{eqnarray*}
p_\lambda(x)&=&\mathbf{P}\bigl(\mbox{$X(t)$ is finite for all $t>0$
and}\\
 &&\hspace*{13pt}\hspace*{6pt}\mbox{does
not go to $\infty$ as $t\to\infty$} \bigr)
\end{eqnarray*}
with $X(0)=x$, as defined in (\ref{eq:gapX}). Note that a time shift
of equation
\eqref{eq:logtan} only changes $\lambda$ and the initial
condition. With
%
\begin{equation}
T=T_\lambda=\frac{4}{\beta} \log\lambda\label{T}
\end{equation}
the
diffusion $\tau\mapsto X(\tau+T)$ satisfies
\eqref{eq:logtan} with $\lambda=1$ and with initial condition
$-\infty$ at $\tau=-T$. This suggests that the dependence on
$\lambda$ for the probability on the right-hand side of
\eqref{eq:gapX} comes mainly from the interval $[0,T]$.
Because of this we take conditional expectations in
\eqref{eq:gapX} with respect to the $\sigma$-algebra
generated by $(X(t), t\in[0,T])$. Using the Markov property
of $X$, we obtain
%
\begin{eqnarray}
p_\lambda=\mathbf{E}\bigl( \mathbf{1}\{ \mbox{$X(t)$ is
finite for all
$0<t\le T$}\} \cdot p_1(X(T))\bigr).
\end{eqnarray}
The~first term in the expectation is a function of the path
$X(t)$ on the time interval $[0,T]$.
Consider a diffusion $Y$ given by the SDE (\ref{eq:Y}) with a drift
function $h(t,x)$ given by Lemma \ref{l_Y}. With the notation of Lemma
\ref{l_Y}, we
set
%
\begin{equation}
\qquad \psi(Y)=\biggl(2-\frac\beta2 \biggr)Y(T)^+ +\frac{\beta}{8}
e^{Y(T)} +\omega(Y(T))+\int_{0}^T \phi\bigl(T-t,Y(t)\bigr)\,dt.\label{eq:psi}
\end{equation}
We apply the
Girsanov transformation of Proposition
\ref{pgirsanov} together with equation (\ref{eq:Girsan}) of
Lemma \ref{l_Y} to get
\[
p_\lambda=\lambda^{\gamma_\beta}e^{-
({\beta}/{64})\lambda^2+({\beta}/{8}-1/4
)\lambda}  \mathbf{E}[ p_1(Y(T_\lambda))
\exp\{\psi(Y)\}],
\]
where $\gamma_\beta=\frac{1}{4}(\frac{\beta}{2}+\frac
{2}{\beta}-
3).$ In order to prove the theorem,
it suffices to show that the
limit
%
\begin{eqnarray}\lim_{\lambda\to\infty}  \mathbf{E}[
p_1(Y(T_\lambda))
\exp\{\psi(Y)\}]\label{eq:lim}
\end{eqnarray}
exists, and is finite and positive. This limit then equals
the constant $\kappa_\beta$ of the asymptotics.
Recall that in (\ref{eq:psi}) the function $\omega$ is continuous and
bounded and $\phi(t,y)$ can be dominated by a function $\tilde\phi
(y)$ which has a finite integral in $[0,\infty)$. 

We will run the process $Y_\lambda(t)$ with a shifted time,
$\tau=t-T=t-\frac4\beta\log\lambda$; that is, let
\[
\tilde Y_T(\tau):=Y_\lambda(\tau+T).
\]
The~advantage of this
shifted time is that the diffusions $\tilde Y_T(\tau)$ for
different $\lambda$ satisfy the same SDE except they evolve
on nested time intervals:
%
\begin{eqnarray}\label{Y()}
d\tilde Y_T(\tau)=\tilde h(\tau,\tilde Y) \,  d\tau+dB,\qquad
\tau>-T,\qquad  \tilde Y_T(-T)=-\infty,
\end{eqnarray}
where the drift term is given by
%
\begin{equation}\label{eq:htilde}
\tilde h(\tau,y)=h(T+\tau,y)=
-\frac{\beta}{8} e^{-\beta\tau/4}
\sinh(y)+h_0(T+\tau,y).
\end{equation}
In this new time-frame,
we need to show that the limit
%
\begin{eqnarray}
\lim_{T \to\infty}  \mathbf{E} p_1(\tilde Y_T(0))
\exp\{\tilde\psi(\tilde Y_T )\}\label{eq:lim1}
\end{eqnarray}
exists, is positive and finite, where
%
\begin{equation}\label{psi}
\tilde\psi(\tilde Y)=\biggl(2-\frac\beta2 \biggr)\tilde Y(0)^+
+\frac{\beta}{8} e^{\tilde Y(0)} +\omega(\tilde
Y(0))+\int_{0}^T \phi(t,\tilde Y(-t))\, dt.
\end{equation}
We will drive the diffusions \eqref{Y()} with the same
Brownian motion $B(t)$. Then for $T_1>T_2$ we have
$Y_{T_1}(\tau)>Y_{T_2}(\tau)$ for $\tau\in[T_2,\infty)$ as
this holds for $\tau=-T_2$ and the domination is preserved
by the evolution.

We also consider a nonnegative-valued diffusion $Z(t)$ given by the SDE
\[
dZ=r(Z) \,dt+dB
\]
which is reflected at $0$ and whose drift term is equal to
%
\begin{equation}\label{eq:driftZ}
r(y)=-\frac{\beta}{16}e^{y}+c_1.
\end{equation}
We will use the stationary version of $Z$ to dominate the diffusions
$\tilde Y_T$.

By Lemma \ref{l_Y}, the term $h_0(y,T+\tau)$ in (\ref{eq:htilde}) is
bounded if $-T\le\tau\le0$. Thus, we can choose the constant $c_1$
in (\ref{eq:driftZ}) so
that
%
\begin{equation}r(z)\ge\sup_{\{\tau<0,  0\le y \le z\}} h(\tau
,y).\label{drift}
\end{equation}
Since $Z$ and $\tilde Y$ are driven by the same Brownian motion, if
$Z,\tilde Y> 0$ then $Z-\tilde Y$ evolves according to
\[
d(Z-\tilde Y)=[r(Z)-f(t,Y)]\, dt.
\]
By (\ref{drift}), this means that if $Z(\tau_0)\ge\tilde Y(\tau_0)$
for a $\tau_0<0$ then this ordering is preserved by the coupling until
time 0.

Consider the process $Z$ in its stationary distribution.
Then $Z(-T)>\break  \tilde Y_T(-T)=-\infty$ therefore $Z$ dominates $\tilde
Y_T$ on $[-T,0]$.
For every fixed $\tau\le0$, the random variables $\tilde
Y_T(\tau)$ are increasing in $T$ and bounded by $Z(\tau)$
so
\[
\tilde Y_\infty(\tau)=\lim_{T\to\infty} \tilde
Y_T(\tau)
\]
exists and is dominated by $Z(\tau)$. The~function
$p_1(x)$ is continuous so $p_1(\tilde Y_T)\to p_1(\tilde Y_{\infty})$.
By \eqref{psi},
we have
\[
\tilde\psi(\tilde Y_T)=a(\tilde Y_T(0))+\int_{0}^T \phi(t,\tilde
Y(-t)) \,dt,
\]
where $a$ is continuous and $\phi(t,y)$ can be dominated by a function
$\tilde\phi(y)$ which has a finite integral in $[0,\infty)$. Hence,
$\tilde\psi(\tilde Y_T)\to\tilde\psi(\tilde Y_\infty)$ and
\[
q_T=e^{\tilde\psi(\tilde Y_T)} p_1(\tilde Y_T)\to
q_\infty=e^{\tilde\psi(\tilde Y_\infty)} p_1(\tilde
Y_\infty)\qquad  \mbox{as }   T\to\infty.
\]
Using Lemma \ref{l:hapci} to estimate $p_1(y)$, we get
\begin{eqnarray*}
q_T \le  c \exp\biggl\{(2-\beta/ 2 )\tilde
Y_T(0)^+
+\frac{\beta}{8} e^{\tilde Y_T(0)}-\frac{\beta}{60} e^{\tilde
Y_T(0)}\biggr\}
  \le c'\chi(\tilde Y_T(0)),
\end{eqnarray*}
where $\chi(y)= \exp\{(\frac{\beta}{8}
-\frac{\beta}{61}) e^{y}\}$.
If we prove that $\mathbf{E}\chi(
Z(0))<\infty$, then the dominated convergence theorem will
imply
%
\begin{equation}\label{szendi}
\mathbf{E}q_T\to\mathbf{E}q_\infty<\infty,
\end{equation}
and the existence of the limiting constant $\kappa_\beta$ will be established.

The~generator of the reflected diffusion $Z$ is given by
\[
\mathcal L f=\tfrac12 f''+f' r
\]
for functions $f$ defined on $[0,\infty)$ with $f'(0+)=0$ [\citet
{RevuzYor}, Chapter VII, Section 3]. Partial integration shows that if
$(\log g)' = 2r $ and $f'(0+)=0$ then
$
\int_0^\infty\mathcal L f(x)   g(x)\, dx=0
$
which means that
\[
g(z)=c\exp(-\beta/8e^z+2c_1z)
\]
gives a stationary density.
Since $\int_0^\infty\chi(z)g(z)\,   dz= \mathbf{E}
\chi(Z(0))<\infty$, the convergence \eqref{szendi} follows.
This shows that
\[
\kappa_\beta= \mathbf{E}q _\infty= \mathbf{E}\bigl[ e^{\psi
(\tilde
Y_\infty)} p_1(\tilde Y_\infty(0))\bigr]<\infty.
\]
The~only thing left to prove is that
$\kappa_\beta= \mathbf{E}  q_\infty$ is not zero. The~definitions
of $q$
and~$\psi$
yield
\[
q_\infty\ge c  p_1(\tilde Y_\infty(0)) e^{(2-\beta/2)
Y_\infty(0)^+}.
\]
By Lemma
\ref{l:hapci}, the function $p_1(\cdot)$ is positive.
Since $\tilde Y_\infty(0)$ is a.s.~finite and we get
that $\mathbf{E}q_\infty>0$ which completes the proof of
Theorem \ref{mainthm}.
\end{pf*}

\section*{Acknowledgments}

We are grateful to Peter Forrester
for introducing us to the original work of \citet{Dy62},
which was the starting point of this paper. We
thank Laure Dumaz, Mu Cai and Yang Chen for helpful comments on a previous
version.

%

\printaddresses

\end{document}